\documentclass[10pt]{amsart}

\usepackage{amssymb,amsmath,amsthm,amsfonts}

\theoremstyle{plain}
\newtheorem{Theorem}{Theorem}

\theoremstyle{definition}
\newtheorem{Remark}[Theorem]{Remark}

\begin{document}
\title[Three remarkable amalgams]
{Three amalgams with remarkable normal subgroup structures}
\author{Diego Rattaggi}
\thanks{Supported by the Swiss National Science Foundation, No.\ PP002--68627}
\address{Universit\'e de Gen\`eve,
Section de math\'ematiques,
2--4 rue du Li\`evre, CP 64, 
CH--1211 Gen\`eve 4, Switzerland}
\email{rattaggi@math.unige.ch}
\date{\today}
\begin{abstract}
We construct three groups $\Lambda_1$, $\Lambda_2$, $\Lambda_3$, which can all be 
decomposed as amalgamated products $F_9 \ast_{F_{81}} F_{9}$
and have very few normal subgroups of finite or infinite index. 
Concretely, $\Lambda_1$ is a simple group, 
$\Lambda_2$ is not simple but has no non-trivial normal subgroup of infinite index, 
and $\Lambda_3$ is not simple but has no proper subgroup of finite index.
\end{abstract}
\maketitle

\section{Introduction} \label{Intro}
Motivated by expected analogies between cocompact lattices 
in products of automorphism groups of regular trees
and cocompact lattices in higher rank semisimple Lie groups,
Burger and Mozes discovered in their study of groups acting on products of trees the first 
examples of finitely presented torsion-free simple groups \cite{BM,BMII}.
These groups are moreover amalgamated products of finitely generated non-abelian free groups, thus
answering Neumann's question \cite{Neumann} on the existence of simple amalgams of free groups. 
One crucial step in the construction of Burger-Mozes is a deep theorem,
which states that certain cocompact lattices in 
the product of automorphism groups of locally finite trees 
$\mathrm{Aut}(T_{1}) \times \mathrm{Aut}(T_{2})$ cannot have non-trivial 
normal subgroups of infinite index.
Applying this theorem to a cocompact lattice which contains as a subgroup a non-residually finite group
constructed by Wise in \cite{Wise}, we give an example of a finitely presented torsion-free simple group
$\Lambda_1$ of the form $F_9 \ast_{F_{81}} F_{9}$, where $F_k$ denotes the free group of rank $k$.
See \cite{RatII} for a list of $32$ other finitely presented torsion-free simple groups 
emerging from the same method.
Note that the simple groups of Burger-Mozes are also explicitly given in principle,
but not very manageable in practice, because of their extremely long finite presentations.
In addition to the simple group $\Lambda_1$, 
we construct two other groups $\Lambda_2$ and $\Lambda_3$, also having amalgam decompositions 
$F_9 \ast_{F_{81}} F_{9}$. 
They are not simple, but $\Lambda_2$ is virtually simple and $\Lambda_3$ has no non-trivial
finite quotients. An amalgam $F_3 \ast_{F_{13}} F_{3}$ without proper subgroups of finite index
has already been constructed by Bhattacharjee in \cite{Bhatta}, using different techniques.
Our search for groups with the desired properties was made 
possible by several computer programs written in GAP \cite{GAP}.
See \cite[Appendix~B]{RatI} for the program code used to construct the examples.
We refer to \cite{BMI}, \cite{BMII}, \cite{RatI} and \cite{Wise}
for detailed background on automorphism groups of trees, lattices in products of trees, 
and square complexes.

\section{Definition of the groups $\Gamma_i$ and $\Lambda_i$} \label{Def}
Let always $i \in \{ 1, 2, 3\}$.
Our groups $\Lambda_i$ will be normal subgroups of index $4$ of groups $\Gamma_i$ defined by their
finite presentations
\[
\Gamma_i = \langle a_1, a_2, a_3, a_4, a_5, b_1, b_2, b_3, b_4, b_5 \mid r_1, \ldots , r_{25} \rangle,
\]
where the relators $r_1, \ldots , r_{25}$ (depending on $i$) are given in Table \ref{Table1}. 
Capital letters in this table indicate inverses, for example 
$r_1 = a_1 b_1 A_2 B_2 = a_1 b_1 a_2^{-1} b_2^{-1}$. 

\begin{table}[h!]
\[
\begin{tabular}{| l | l || l | l | l | l |}
\hline
         & $\Gamma_1$, $\Gamma_2$, $\Gamma_3$  &  & $\Gamma_1$ &  $\Gamma_2$ &  $\Gamma_3$  \\ \hline
$r_1$    & $a_1 b_1 A_2 B_2$ & 
$r_{13}$ & $a_1 b_4 a_2 B_5$  &  $a_1 b_4 A_4 b_5$  & $a_1 b_4 a_1 b_5$ \\ \hline
$r_2$    & $a_1 b_2 A_1 B_1$ &
$r_{14}$ & $a_1 b_5 A_5 b_4$  &  $a_1 b_5 a_2 B_5$  & $a_1 B_5 a_2 B_5$ \\ \hline
$r_3$    & $a_1 b_3 A_2 B_3$ &
$r_{15}$ & $a_1 B_5 a_3 B_4$  &  $a_1 B_5 a_3 B_4$  & $a_1 B_4 A_4 B_4$ \\ \hline
$r_4$    & $a_1 B_3 A_2 b_2$ & 
$r_{16}$ & $a_1 B_4 a_3 b_5$  &  $a_1 B_4 a_2 b_4$  & $a_2 b_4 a_2 b_5$ \\ \hline
$r_5$    & $a_1 B_1 A_2 b_3$ &
$r_{17}$ & $a_2 b_4 A_2 b_5$  &  $a_2 b_5 A_3 b_4$  & $a_2 B_4 A_3 B_4$ \\ \hline
$r_6$    & $a_2 b_2 A_2 B_1$ &
$r_{18}$ & $a_2 b_5 a_4 B_4$  &  $a_2 B_4 a_4 B_5$  & $a_3 b_5 a_4 B_4$ \\ \hline
$r_7$    & $a_3 b_1 A_4 B_2$ & 
$r_{19}$ & $a_3 b_4 a_4 b_5$  &  $a_3 b_4 a_4 b_4$  & $a_3 B_5 A_5 B_5$ \\ \hline
$r_8$    & $a_3 b_2 A_3 B_1$ & 
$r_{20}$ & $a_3 B_5 a_4 b_4$  &  $a_3 b_5 A_5 b_5$  & $a_3 B_4 a_4 b_5$ \\ \hline
$r_9$    & $a_3 b_3 A_4 B_3$ & 
$r_{21}$ & $a_4 B_5 A_5 B_4$  &  $a_4 b_5 a_5 b_5$  & $a_4 B_5 a_5 B_5$ \\ \hline
$r_{10}$ & $a_3 B_3 A_4 b_2$ & 
$r_{22}$ & $a_5 b_1 A_5 b_3$  &  $a_5 b_1 A_5 b_3$  & $a_5 b_1 a_5 b_4$ \\ \hline
$r_{11}$ & $a_3 B_1 A_4 b_3$ &
$r_{23}$ & $a_5 b_2 A_5 B_5$  &  $a_5 b_2 A_5 B_1$  & $a_5 b_2 A_5 b_3$ \\ \hline
$r_{12}$ & $a_4 b_2 A_4 B_1$ &
$r_{24}$ & $a_5 b_3 A_5 B_1$  &  $a_5 b_3 A_5 B_4$  & $a_5 b_3 A_5 b_2$ \\ \hline
         &                   &
$r_{25}$ & $a_5 b_4 A_5 B_2$  &  $a_5 b_4 A_5 B_2$  & $a_5 B_4 a_5 B_1$ \\ \hline
\end{tabular}
\]
\caption{The 25 relators of $\Gamma_1$, $\Gamma_2$, $\Gamma_3$} \label{Table1}
\end{table}

Observe that the twelve relators $r_1, \ldots ,r_{12}$ are the same for each group $\Gamma_i$.
The reason for this will become clear in the proof of Theorem~\ref{Thm} in Section~\ref{Res}.
To describe the geometric nature of $\Gamma_i$, we recall 
the following general construction which associates 
to a finite presentation of a group $G$ its \emph{standard $2$-complex} $X$ with fundamental group $G$:
by definition, the one-skeleton of $X$ has a single vertex $x$ and an oriented loop 
for each generator of the given presentation of $G$.
Furthermore, for each relator $r$, a $2$-cell with boundary labelled by $r$ 
is glued into this one-skeleton to get $X$. Then $G = \pi_1(X,x)$.
By construction of the $25$ relators of $\Gamma_i$, 
its associated standard $2$-complex $X_i$ is a finite square complex
(all relators have length four, hence all $2$-cells are squares)
having the additional property that its universal cover $\tilde{X_i}$ 
is the affine building $\mathcal{T}_{10} \times \mathcal{T}_{10}$,
the product of two $10$-regular trees.
Equivalently, this property requires that to each pair $(a,b) \in A \times B$,
there is a uniquely determined pair $(\tilde{a}, \tilde{b}) \in A \times B$ such that
$ab = \tilde{b} \tilde{a}$ in $\Gamma_i$,
where $A := \{ a_1, \ldots , a_5 \}^{\pm 1}$ and $B := \{ b_1, \ldots , b_5 \}^{\pm 1}$.
This can be easily verified for our three given examples.
In the terminology of \cite{BMII}, $X_i$ is a finite $1$-vertex VH-T-square complex,
and in the terminology of \cite{RatI,RatII}, $\Gamma_i = \pi_1(X_i)$ is a $(10,10)$--group.
The group of automorphisms $\mathrm{Aut}(\mathcal{T}_{10})$, 
equipped with the usual topology of simple convergence,
is a locally compact group. Taking the product topology, $\Gamma_i$ can be seen as a discrete subgroup of
$\mathrm{Aut}(\mathcal{T}_{10}) \times \mathrm{Aut}(\mathcal{T}_{10})$ with compact quotient, 
in other words as a cocompact lattice.
A crucial role in deducing interesting results on the normal subgroup structure of $\Gamma_i$
play the so-called local groups of $\Gamma_i$. The idea to define them is the following:
take the projection of $\Gamma_i$ to one factor of 
$\mathrm{Aut}(\mathcal{T}_{10}) \times \mathrm{Aut}(\mathcal{T}_{10})$
(say the projection $\mathrm{pr}_1$ to the first factor) 
and fix any vertex $x_h$ of $\mathcal{T}_{10}$.
Then the elements in the closure $\overline{\mathrm{pr}_1(\Gamma_i)} < \mathrm{Aut}(\mathcal{T}_{10})$
stabilizing $x_h$, induce a finite permutation group $P_h^{(1)}(\Gamma_i) < S_{10}$ 
on the $10$ neighbouring vertices
of $x_h$ in $\mathcal{T}_{10}$ (or more generally, for $k \in \mathbb{N}$, 
subgroups $P_h^{(k)}(\Gamma_i)$ of the symmetric group
$S_{10 \cdot 9^{k-1}}$, taking the induced action on the $k$-sphere in $\mathcal{T}_{10}$
around $x_h$). 
The same procedure can be done with the second projection $\mathrm{pr}_2$ to get local groups 
$P_v^{(k)}(\Gamma_i) < S_{10 \cdot 9^{k-1}}$.
It is important to note that these local groups 
(more precisely, their generators in $S_{10 \cdot 9^{k-1}}$)
can be directly computed, given the relators $r_1, \ldots, r_{25}$ of Table~\ref{Table1}, 
see \cite[Chapter~1]{BMII} or \cite[Section~1.4]{RatI} for details.
Here, we get for $k = 1$ the groups
\begin{align}
P_h^{(1)}(\Gamma_1) = \langle &(7,8)(9,10), (1,2)(3,4), (1,2)(3,4)(7,8)(9,10), \notag \\
                              &(1,8,4,5)(2,7,3,10), (1,9,4,8)(3,10,6,7) \rangle = A_{10}, \notag
\end{align}
\begin{align}
P_h^{(1)}(\Gamma_2) = \langle &(7,8)(9,10), (1,2)(3,4), (1,2)(3,4)(7,8)(9,10), \notag \\
                              &(1,8,4,9)(2,10,7,3), (1,9,8,6,4)(2,7,5,3,10) \rangle = A_{10}, \notag
\end{align}
\begin{align}
P_h^{(1)}(\Gamma_3) = \langle &(5,6)(7,8)(9,10), (1,2)(3,4), (1,2)(3,4)(7,8)(9,10), \notag \\
                              &(1,4,8,9,2,3,7,10)(5,6), (1,9,2,10)(3,5,7)(4,6,8) \rangle, \notag
\end{align}
\begin{align}
P_v^{(1)}(\Gamma_1) = \langle &(1,2)(4,6,7,5)(8,10,9), (1,2,3)(4,5,7,6)(9,10), (1,2)(4,5,7,6)(8,10,9), \notag \\
                              &(1,2,3)(4,6,7,5)(9,10), (1,3,10,8)(2,4,6,9,7,5) \rangle = A_{10}, \notag
\end{align}
\begin{align}
P_v^{(1)}(\Gamma_2) = \langle &(1,2)(4,6)(8,10,9), (1,2,3)(5,7)(9,10), (1,2)(4,6,5,7)(8,10,9), \notag \\
                              &(1,2,3)(4,6,5,7)(9,10), (1,2,4,3,10,9,7,8)(5,6) \rangle = A_{10}, \notag
\end{align}
\begin{align}
P_v^{(1)}(\Gamma_3) = \langle &(1,2)(4,7,5,6)(8,10,9), (1,2,3)(4,7,5,6)(9,10), (1,2)(4,5,6,7)(8,10,9), \notag \\
                              &(1,2,3)(4,5,6,7)(9,10), (1,7)(2,8)(3,9)(4,10)(5,6) \rangle = S_{10}. \notag
\end{align}
The transitivity of the permutation groups given above will be important
in the proof of Theorem~\ref{Thm}. Recall that a group $G < S_{10}$ is
\emph{transitive} if for any pair $m, n \in \{1, \ldots, 10 \}$ there exists
a $g \in G$ such that $g(m) = n$. Moreover,
$G$ is called \emph{$2$-transitive} if for any $m_1, m_2, n_1, n_2 \in \{1, \ldots, 10 \}$
with $m_1 \ne m_2$ and $n_1 \ne n_2$ there is an element $g \in G$ such that
$g(m_1) = n_1$ and $g(m_2) = n_2$.
Note that the group $P_h^{(1)}(\Gamma_3)$ is a transitive 
(but not $2$-transitive) subgroup of $S_{10}$ of order $3840$,
whereas the alternating group $A_{10}$ and the symmetric group $S_{10}$
are obviously $2$-transitive.

We define now $\Lambda_i$ to be the kernel of the surjective homomorphism
\begin{align}
\Gamma_i &\to \mathbb{Z} / 2\mathbb{Z} \times \mathbb{Z} / 2\mathbb{Z} \notag \\
a_1, \ldots , a_5 &\mapsto (1 + 2\mathbb{Z}, 0 + 2\mathbb{Z}) \notag \\
b_1, \ldots , b_5 &\mapsto (0 + 2\mathbb{Z}, 1 + 2\mathbb{Z}), \notag
\end{align}
where $\Gamma_i$ is given by its finite presentation described above.
Each group $\Lambda_i$ can be decomposed in two ways as amalgamated products
$F_9 \ast_{F_{81}} F_{9}$, such that $F_{81}$ has index $10$ in both factors $F_9$.
More precisely, this means that for any $i \in \{1,2,3\}$ there exist
injective homomorphisms $j_1, j_3 : F_{81} \to F_9 \cong \langle s_1, \ldots, s_9 \rangle$
and $j_2, j_4 : F_{81} \to F_9 \cong \langle t_1, \ldots, t_9 \rangle$
such that 
\[
[F_9 : j_1 (F_{81}) ] =  [F_9 : j_2 (F_{81}) ] =  [F_9 : j_3 (F_{81}) ] =  [F_9 : j_4 (F_{81}) ] = 10
\]
and
\begin{align}
\Lambda_i &\cong \langle s_1, \ldots, s_9, t_1, \ldots, t_9 \mid 
j_1(u_1) = j_2(u_1), \ldots, j_1(u_{81}) = j_2(u_{81}) \rangle \notag \\
&\cong \langle s_1, \ldots, s_9, t_1, \ldots, t_9 \mid 
j_3(u_1) = j_4(u_1), \ldots, j_3(u_{81}) = j_4(u_{81}) \rangle, \notag
\end{align}
where $\{u_1, \ldots, u_{81}\}$ are the free generators of $F_{81}$.
This is a direct consequence of a result of Wise (see \cite[Theorem I.1.18]{Wise}),
describing each of the two decompositions of certain square complex groups $\Gamma$ as a 
fundamental group of a finite graph of finitely generated free groups
(in the language of the Bass-Serre theory).
If the local groups of $\Gamma$ are ``sufficiently transitive'' (which always happens in our examples),
the two finite graphs corresponding to $\Lambda_i$ in Wise's construction each consist of 
two vertices and a single edge.
Therefore we get amalgams of finitely generated free groups.
It is well-known that amalgams of free groups are always torsion-free, 
since every element of finite order in an amalgam is conjugate to an element 
of finite order in one of the two factors (see for example \cite[Theorem~IV.2.7]{LS}).
Note that following Wise's proof of \cite[Theorem I.1.18]{Wise}, it is not difficult
(but quite laborious by hand) to give explicit descriptions of the 
injective homomorphisms $F_{81} \to F_9$ in the amalgam decompositions of $\Lambda_i$. 

\section{Results and Proofs} \label{Res}
In the following theorem, we discuss the normal subgroups of $\Lambda_i$.
\begin{Theorem} \label{Thm}
Let $\Lambda_1$, $\Lambda_2$, $\Lambda_3$ be the groups defined in Section~\ref{Def}. Then
\begin{itemize}
\item[$(1)$] $\Lambda_1$ is simple.
\item[$(2)$] Every non-trivial normal subgroup of $\Lambda_2$ has finite index, but $\Lambda_2$ is not simple.
\item[$(3)$] $\Lambda_3$ has no proper subgroups of finite index, but is not simple. 
\end{itemize}
\end{Theorem}
\begin{proof}
Let $W$ be the group with finite presentation 
\[
\langle a_1, a_2, a_3, a_4, b_1, b_2, b_3 \mid r_1, \ldots , r_{12} \rangle,
\]
where the relators $r_1, \ldots , r_{12}$ are again taken from Table \ref{Table1}.
Wise showed in \cite[Main Theorem 5.5]{Wise}, that the non-trivial element 
$w := a_2 a_1^{-1} a_3 a_4^{-1} \in W$ is contained in each
finite index subgroup of $W$. In particular, $W$ is non-residually finite.
Moreover, $W < \mathrm{Aut}(\mathcal{T}_{8}) \times \mathrm{Aut}(\mathcal{T}_{6})$ 
is the fundamental group of a $1$-vertex VH-T square complex which embeds
into the square complex $X_i$ associated to $\Gamma_i$ ($i = 1,2,3$), inducing an injection
on the level of fundamental groups, i.e.\ $W < \Gamma_i = \pi_1(X_i)$ 
(the fact that we get an injection can be deduced from the non-positive curvature 
of the product of trees $\mathcal{T}_{10} \times \mathcal{T}_{10}$, 
see \cite[Proposition~II.4.14(1)]{BrHa}). 
Hence we have
\[
1 \ne w \in  \bigcap_{N \overset{\text{f.i.}}{\lhd} W} N < 
\bigcap_{N \overset{\text{f.i.}}{\lhd} \Gamma_i} N = 
\bigcap_{N \overset{\text{f.i.}}{<} \Gamma_i} N \: \lhd \: \Gamma_i,
\]
where ``f.i.'' stands for ``finite index''.
In particular, $\Gamma_i$ (and hence its finite index subgroup $\Lambda_i$) is
non-residually finite.
Observe that $w \in \Lambda_i \lhd \Gamma_i$.
One important point in the construction of $\Gamma_i$ is to guarantee 
that the normal closure of $w$ in $\Gamma_i$, 
denoted by
$\langle \! \langle w \rangle \! \rangle_{\Gamma_i}$, has finite index in $\Lambda_i$.
(Note that however $[W : \langle \! \langle w \rangle \! \rangle_W] = \infty$.) 
This already implies that
$\langle \! \langle w \rangle \! \rangle_{\Gamma_i}$ has no proper subgroups of finite index.
Indeed, assume that $M < \langle \! \langle w \rangle \! \rangle_{\Gamma_i}$ 
is a subgroup of finite index. Then
\[
\bigcap_{N \overset{\text{f.i.}}{<} \Gamma_i} N < M \overset{\text{f.i.}}{<} 
\langle \! \langle w \rangle \! \rangle_{\Gamma_i}
\overset{\text{f.i.}}{<} \Lambda_i \overset{\text{f.i.}}{<} \Gamma_i .
\]
Using 
\[
\langle \! \langle w \rangle \! \rangle_{\Gamma_i} 
< \bigcap_{N \overset{\text{f.i.}}{\lhd} \Gamma_i} N
= \bigcap_{N \overset{\text{f.i.}}{<} \Gamma_i} N,
\]
we get
\[
M = \langle \! \langle w \rangle \! \rangle_{\Gamma_i} = \bigcap_{N \overset{\text{f.i.}}{<} \Gamma_i} N.
\]
We proceed now separately for the three groups $\Lambda_1$, $\Lambda_2$ and $\Lambda_3$.
\begin{itemize}
\item[$(1)$]We have $\langle \! \langle w \rangle \! \rangle_{\Gamma_1} = \Lambda_1$. 
This can be checked by hand, or more easily,
using a computer algebra system like GAP \cite{GAP}, 
which shows that adding the relator $w$ to the presentation of
$\Gamma_1$ gives the group $\mathbb{Z} / 2\mathbb{Z} \times \mathbb{Z} / 2\mathbb{Z}$ of order $4$. 
It remains to prove that $\Lambda_1$ has no non-trivial normal subgroups of 
infinite index. But this follows directly from the normal subgroup theorem of 
Burger-Mozes \cite[Theorem 4.1, Corollary 5.4]{BMII}
applied to the ``irreducible'' cocompact lattice 
$\Gamma_1 < \mathrm{Aut}(\mathcal{T}_{10}) \times \mathrm{Aut}(\mathcal{T}_{10})$
with local groups $P_h^{(1)}(\Gamma_1) \cong P_v^{(1)}(\Gamma_1) = A_{10}$, and applied to 
its finite index subgroup $\Lambda_1 < \Gamma_1$.
\item[$(2)$]For the second group, we compute 
$[\Lambda_2 : \langle \! \langle w \rangle \! \rangle_{\Gamma_2}] = 2$, 
thus $\Lambda_2$ is not simple.
By exactly the same argument as in part~(1), every non-trivial normal subgroup of $\Gamma_2$ 
(and of $\Lambda_2$, respectively) has finite index. Observe that
$\langle \! \langle w \rangle \! \rangle_{\Gamma_2}$ is a simple group with amalgam 
decomposition $F_{17} \ast_{F_{161}} F_{17}$.
In particular, $\Gamma_2$ and $\Lambda_2$ are virtually simple groups.
\item[$(3)$]As in part~(1), $\langle \! \langle w \rangle \! \rangle_{\Gamma_3} = \Lambda_3$ proves 
that $\Lambda_3$ has no proper subgroup of finite index.
However, in contrast to what happens in part~(1) and (2), the local group $P_h^{(1)}(\Gamma_3)$
is transitive, but not $2$-transitive. Therefore, the normal subgroup theorem of Burger-Mozes
cannot be applied here. Indeed, $\Lambda_3$ is not simple, since
$1 \ne \langle \! \langle a_5^4 \rangle \! \rangle_{\Lambda_3} \ne \Lambda_3$. 
This comes from the fact that $a_5^4$ acts trivially
on the second factor of $\mathcal{T}_{10} \times \mathcal{T}_{10}$. In other words, 
$a_5^4 \in \mathrm{ker}(\mathrm{pr}_2) \lhd \Gamma_3$.
To see this, let 
\[
A' := \{ (a_1 a_2^{-1})^{2}, (a_2^{-1} a_1)^{2}, (a_3 a_4^{-1})^{2}, (a_4^{-1} a_3)^{2}, a_5 ^{4} \}^{\pm 1}
\]
and check that for all $a' \in A'$ and $b \in B = \{ b_1, \ldots, b_5 \}^{\pm 1}$, we have
$b^{-1} a' b \in A'$. This in fact implies that $A' \subset \mathrm{ker}(\mathrm{pr}_2)$.
Note that no element of $\Gamma_3$ acts trivially on the \emph{first} factor of
 $\mathcal{T}_{10} \times \mathcal{T}_{10}$
(by \cite[Proposition~3.1.2,~1)]{BMI} and \cite[Proposition~3.3.2]{BMI}).
As a consequence, $\Lambda_3$ has two decompositions $F_9 \ast_{F_{81}} F_{9}$, 
where one amalgam is effective and the other one is not effective. 
\end{itemize}
\end{proof}

We conclude by giving two remarks:
\begin{Remark}
Recall that a group $G$ is called \emph{SQ-universal} 
if every countable group can be embedded in a quotient of $G$.
It is mentioned in \cite[Chapter 9.15]{BassLub} that Ilya Rips can prove any amalgamated
product $A \ast_{C} B$ to be SQ-universal, 
provided that $B \ne C$ and the number of double cosets $|C \backslash A / C|$
is at least $3$ (if $C$ is seen as usual as a subgroup of $A$ and $B$ via the two injections 
$j_1: C \to A$ and $j_2: C \to B$ in the amalgam),
but there is no published proof as far as we know. 
If Rips' statement is true, we could apply it to exactly 
one decomposition $F_9 \ast_{F_{81}} F_{9}$ of $\Lambda_3$ (to the effective one), where
$|F_{81} \backslash F_9 / F_{81}| = 3$. Note however that in the second decomposition of $\Lambda_3$ 
(where the corresponding local group $P_v^{(1)}(\Gamma_3)$ is $S_{10}$) and in 
both decompositions of $\Lambda_1$ and $\Lambda_2$, we always have $|F_{81} \backslash F_9 / F_{81}| = 2$, 
since their local actions on $\mathcal{T}_{10}$ are $2$-transitive.
\end{Remark}

\begin{Remark}
By construction, the three groups $\Lambda_1$, $\Lambda_2$, $\Lambda_3$ are non-residually finite.
As a contrast, if one takes a \emph{double} $F_9 \ast_{F_{81}} F_{9}$ 
(i.e.\ an amalgam where the two injections $j_1, j_2 : F_{81} \to F_9$ are identical),
such that $F_{81}$ has finite index in both factors $F_9$ 
(consequently index $10 = (81-1)/(9-1)$),
then one directly gets a surjective homomorphism $F_9 \ast_{F_{81}} F_{9} \to F_9$
(the obvious folding map),
and moreover $F_9 \ast_{F_{81}} F_{9}$ contains by \cite[Theorem~1.4]{BDGM}
a subgroup of finite index which is a direct product of two non-abelian free groups of finite rank.
In particular, such a double $F_9 \ast_{F_{81}} F_{9}$ is SQ-universal and residually finite.
The residual finiteness also follows from \cite{Stebe}.
\end{Remark}

\end{document}